\documentclass[12pt]{amsart}
 \usepackage{amsmath}
\usepackage{amssymb,latexsym}
\usepackage{amscd}
\usepackage{comment}
 \usepackage{xspace}
 \usepackage{verbatim}
\usepackage[english]{babel}

\begin{document}

\numberwithin{equation}{section}
\title[Semilinear elliptic equations with Hardy potential]{Moderate
solutions  of semilinear elliptic equations with Hardy potential}
\author{Moshe Marcus }
\address{Department of Mathematics, Technion\\
 Haifa 32000, ISRAEL}
 \email{marcusm@math.technion.ac.il}
\author{Phuoc-Tai Nguyen}
\address{Department of Mathematics, Technion\\
 Haifa 32000, ISRAEL}
 \email{nguyenphuoctai.hcmup@gmail.com}

\date{}

\maketitle


\newcommand{\txt}[1]{\;\text{ #1 }\;}
\newcommand{\tbf}{\textbf}
\newcommand{\tit}{\textit}
\newcommand{\tsc}{\textsc}
\newcommand{\trm}{\textrm}
\newcommand{\mbf}{\mathbf}
\newcommand{\mrm}{\mathrm}
\newcommand{\bsym}{\boldsymbol}
\newcommand{\scs}{\scriptstyle}
\newcommand{\sss}{\scriptscriptstyle}
\newcommand{\txts}{\textstyle}
\newcommand{\dsps}{\displaystyle}
\newcommand{\fnz}{\footnotesize}
\newcommand{\scz}{\scriptsize}
\newcommand{\be}{\begin{equation}}
\newcommand{\bel}[1]{\begin{equation}\label{#1}}
\newcommand{\ee}{\end{equation}}
\newcommand{\eqnl}[2]{\begin{equation}\label{#1}{#2}\end{equation}}
\newcommand{\barr}{\begin{eqnarray}}
\newcommand{\earr}{\end{eqnarray}}
\newcommand{\bars}{\begin{eqnarray*}}
\newcommand{\ears}{\end{eqnarray*}}
\newcommand{\nnu}{\nonumber \\}

\newtheorem{subn}{\name}
\renewcommand{\thesubn}{}
\newcommand{\bsn}[1]{\def\name{#1}\begin{subn}}
\newcommand{\esn}{\end{subn}}

\newtheorem{sub}{\name}[section]
\newcommand{\dn}[1]{\def\name{#1}}


\newcommand{\bs}{\begin{sub}}
\newcommand{\es}{\end{sub}}
\newcommand{\bsl}[1]{\begin{sub}\label{#1}}

\newcommand{\bth}[1]{\def\name{Theorem}
\begin{sub}\label{t:#1}}
\newcommand{\blemma}[1]{\def\name{Lemma}
\begin{sub}\label{l:#1}}
\newcommand{\bcor}[1]{\def\name{Corollary}
\begin{sub}\label{c:#1}}
\newcommand{\bdef}[1]{\def\name{Definition}
\begin{sub}\label{d:#1}}
\newcommand{\bprop}[1]{\def\name{Proposition}
\begin{sub}\label{p:#1}}
\newcommand{\R}{\eqref}
\newcommand{\rth}[1]{Theorem~\ref{t:#1}}
\newcommand{\rlemma}[1]{Lemma~\ref{l:#1}}
\newcommand{\rcor}[1]{Corollary~\ref{c:#1}}
\newcommand{\rdef}[1]{Definition~\ref{d:#1}}
\newcommand{\rprop}[1]{Proposition~\ref{p:#1}}
\newcommand{\BA}{\begin{array}}
\newcommand{\EA}{\end{array}}
\newcommand{\BAN}{\renewcommand{\arraystretch}{1.2}
\setlength{\arraycolsep}{2pt}\begin{array}}
\newcommand{\BAV}[2]{\renewcommand{\arraystretch}{#1}
\setlength{\arraycolsep}{#2}\begin{array}}
\newcommand{\BSA}{\begin{subarray}}
\newcommand{\ESA}{\end{subarray}}


\newcommand{\BAL}{\begin{aligned}}
\newcommand{\EAL}{\end{aligned}}
\newcommand{\BALG}{\begin{alignat}}
\newcommand{\EALG}{\end{alignat}}
\newcommand{\BALGN}{\begin{alignat*}}
\newcommand{\EALGN}{\end{alignat*}}
\newcommand{\note}[1]{\noindent\textit{#1.}\hspace{2mm}}
\newcommand{\Proof}{\note{Proof}}
\newcommand{\Remark}{\note{Remark}}
\newcommand{\modin}{$\,$\\[-4mm] \indent}
\newcommand{\forevery}{\quad \forall}
\newcommand{\set}[1]{\{#1\}}
\newcommand{\setdef}[2]{\{\,#1:\,#2\,\}}
\newcommand{\setm}[2]{\{\,#1\mid #2\,\}}
\newcommand{\mt}{\mapsto}
\newcommand{\lra}{\longrightarrow}
\newcommand{\lla}{\longleftarrow}
\newcommand{\llra}{\longleftrightarrow}
\newcommand{\Lra}{\Longrightarrow}
\newcommand{\Lla}{\Longleftarrow}
\newcommand{\Llra}{\Longleftrightarrow}
\newcommand{\warrow}{\rightharpoonup}

\newcommand{\paran}[1]{\left (#1 \right )}
\newcommand{\sqbr}[1]{\left [#1 \right ]}
\newcommand{\curlybr}[1]{\left \{#1 \right \}}
\newcommand{\abs}[1]{\left |#1\right |}

\newcommand{\norm}[1]{\left \|#1\right \|}

\newcommand{\paranb}[1]{\big (#1 \big )}

\newcommand{\lsqbrb}[1]{\big [#1 \big ]}

\newcommand{\lcurlybrb}[1]{\big \{#1 \big \}}

\newcommand{\absb}[1]{\big |#1\big |}

\newcommand{\normb}[1]{\big \|#1\big \|}

\newcommand{
\paranB}[1]{\Big (#1 \Big )}

\newcommand{\absB}[1]{\Big |#1\Big |}

\newcommand{\normB}[1]{\Big \|#1\Big \|}
\newcommand{\produal}[1]{\langle #1 \rangle}

\newcommand{\thkl}{\rule[-.5mm]{.3mm}{3mm}}
\newcommand{\thknorm}[1]{\thkl #1 \thkl\,}
\newcommand{\trinorm}[1]{|\!|\!| #1 |\!|\!|\,}
\newcommand{\bang}[1]{\langle #1 \rangle}
\def\angb<#1>{\langle #1 \rangle}
\newcommand{\vstrut}[1]{\rule{0mm}{#1}}
\newcommand{\rec}[1]{\frac{1}{#1}}
\newcommand{\opname}[1]{\mbox{\rm #1}\,}
\newcommand{\supp}{\opname{supp}}
\newcommand{\dist}{\opname{dist}}
\newcommand{\myfrac}[2]{{\displaystyle \frac{#1}{#2} }}
\newcommand{\myint}[2]{{\displaystyle \int_{#1}^{#2}}}
\newcommand{\mysum}[2]{{\displaystyle \sum_{#1}^{#2}}}
\newcommand {\dint}{{\displaystyle \myint\!\!\myint}}
\newcommand{\q}{\quad}
\newcommand{\qq}{\qquad}
\newcommand{\hsp}[1]{\hspace{#1mm}}
\newcommand{\vsp}[1]{\vspace{#1mm}}
\newcommand{\ity}{\infty}
\newcommand{\prt}{\partial}
\newcommand{\sms}{\setminus}
\newcommand{\ems}{\emptyset}
\newcommand{\ti}{\times}
\newcommand{\pr}{^\prime}
\newcommand{\ppr}{^{\prime\prime}}
\newcommand{\tl}{\tilde}
\newcommand{\sbs}{\subset}
\newcommand{\sbeq}{\subseteq}
\newcommand{\nind}{\noindent}
\newcommand{\ind}{\indent}
\newcommand{\ovl}{\overline}
\newcommand{\unl}{\underline}
\newcommand{\nin}{\not\in}
\newcommand{\pfrac}[2]{\genfrac{(}{)}{}{}{#1}{#2}}


\def\ga{\alpha}     \def\gb{\beta}       \def\gg{\gamma}
\def\gc{\chi}       \def\gd{\delta}      \def\ge{\epsilon}
\def\gth{\theta}                         \def\vge{\varepsilon}
\def\gf{\phi}       \def\vgf{\varphi}    \def\gh{\eta}
\def\gi{\iota}      \def\gk{\kappa}      \def\gl{\lambda}
\def\gm{\mu}        \def\gn{\nu}         \def\gp{\pi}
\def\vgp{\varpi}    \def\gr{\gd}        \def\vgr{\varrho}
\def\gs{\sigma}     \def\vgs{\varsigma}  \def\gt{\tau}
\def\gu{\upsilon}   \def\gv{\vartheta}   \def\gw{\omega}
\def\gx{\xi}        \def\gy{\psi}        \def\gz{\zeta}
\def\Gg{\Gamma}     \def\Gd{\Delta}      \def\Gf{\Phi}
\def\Gth{\Theta}
\def\Gl{\Lambda}    \def\Gs{\Sigma}      \def\Gp{\Pi}
\def\Gw{\Omega}     \def\Gx{\Xi}         \def\Gy{\Psi}

\def\CS{{\mathcal S}}   \def\CM{{\mathcal M}}   \def\CN{{\mathcal N}}
\def\CR{{\mathcal R}}   \def\CO{{\mathcal O}}   \def\CP{{\mathcal P}}
\def\CA{{\mathcal A}}   \def\CB{{\mathcal B}}   \def\CC{{\mathcal C}}
\def\CD{{\mathcal D}}   \def\CE{{\mathcal E}}   \def\CF{{\mathcal F}}
\def\CG{{\mathcal G}}   \def\CH{{\mathcal H}}   \def\CI{{\mathcal I}}
\def\CJ{{\mathcal J}}   \def\CK{{\mathcal K}}   \def\CL{{\mathcal L}}
\def\CT{{\mathcal T}}   \def\CU{{\mathcal U}}   \def\CV{{\mathcal V}}
\def\CZ{{\mathcal Z}}   \def\CX{{\mathcal X}}   \def\CY{{\mathcal Y}}
\def\CW{{\mathcal W}} \def\CQ{{\mathcal Q}}
\def\BBA {\mathbb A}   \def\BBb {\mathbb B}    \def\BBC {\mathbb C}
\def\BBD {\mathbb D}   \def\BBE {\mathbb E}    \def\BBF {\mathbb F}
\def\BBG {\mathbb G}   \def\BBH {\mathbb H}    \def\BBI {\mathbb I}
\def\BBJ {\mathbb J}   \def\BBK {\mathbb K}    \def\BBL {\mathbb L}
\def\BBM {\mathbb M}   \def\BBN {\mathbb N}    \def\BBO {\mathbb O}
\def\BBP {\mathbb P}   \def\BBR {\mathbb R}    \def\BBS {\mathbb S}
\def\BBT {\mathbb T}   \def\BBU {\mathbb U}    \def\BBV {\mathbb V}
\def\BBW {\mathbb W}   \def\BBX {\mathbb X}    \def\BBY {\mathbb Y}
\def\BBZ {\mathbb Z}

\def\GTA {\mathfrak A}   \def\GTB {\mathfrak B}    \def\GTC {\mathfrak C}
\def\GTD {\mathfrak D}   \def\GTE {\mathfrak E}    \def\GTF {\mathfrak F}
\def\GTG {\mathfrak G}   \def\GTH {\mathfrak H}    \def\GTI {\mathfrak I}
\def\GTJ {\mathfrak J}   \def\GTK {\mathfrak K}    \def\GTL {\mathfrak L}
\def\GTM {\mathfrak M}   \def\GTN {\mathfrak N}    \def\GTO {\mathfrak O}
\def\GTP {\mathfrak P}   \def\GTR {\mathfrak R}    \def\GTS {\mathfrak S}
\def\GTT {\mathfrak T}   \def\GTU {\mathfrak U}    \def\GTV {\mathfrak V}
\def\GTW {\mathfrak W}   \def\GTX {\mathfrak X}    \def\GTY {\mathfrak Y}
\def\GTZ {\mathfrak Z}   \def\GTQ {\mathfrak Q}

\def\sign{\mathrm{sign\,}}
\def\bdw{\prt\Gw\xspace}
\def\nabu{|\nabla u|}
\def\tr{\mathrm{tr\,}}
\def\gap{{\ga_+}}
\def\gan{{\ga_-}}
\tableofcontents

\begin{abstract}
Let $\Gw$ be a bounded smooth domain in $\BBR^N$.
We study positive solutions of equation (E) $-L_\mu u+ u^q = 0$ in $\Gw$ where $L_\mu=\Gd + \frac{\gm}{\gd^2}$,
 $0<\gm$, $q>1$ and $\gd(x)=\mathrm{dist}\,(x,\bdw)$. A positive solution of (E) is moderate if it is dominated by an $L_\mu$-harmonic function. If $\mu<C_H(\Gw)$ (the Hardy constant for $\Gw$) every positive $L_\mu$- harmonic functions can be represented in terms of a finite measure on $\bdw$ via the Martin representation theorem. However the classical measure boundary trace  of any such solution is zero. We introduce a notion of  normalized boundary  trace by which we obtain a complete classification of the positive moderate solutions of (E) in the subcritical case, $1<q<q_{\mu,c}$. (The critical value depends only on $N$ and $\mu$.) For $q\geq q_{\mu,c}$ there exists no moderate solution with an isolated singularity on the boundary. The normalized boundary trace and associated boundary value problems are also discussed in detail for the linear operator $L_\mu$. These results form the basis for the study of the nonlinear problem.

\medskip


\end{abstract}


\section{Introduction}

In this paper, we investigate boundary value problem with measure data for the following equations
\bel{N} - \Gd u - \myfrac{\gm}{\gd^2}u + u^q = 0 \ee
in a $C^2$ bounded domain $\Gw$, where $q>1$, $\gm \in \BBR$ and $\gd(x)=\dist(x,\prt \Gw)$. This problem is naturally linked to the theory of linear Schr\"odinger equations $-L^Vu=0$ where
$L^V := \Gd +V$ and the potential $V$ satisfies
 $|V|\leq c\gd^{-2}$.
Such equations have been studied in numerous papers (see e.g. \cite{An1, An2} and the references therein).

Put
\begin{equation}\label{L}
  L_\mu:=\Gd+\frac{\mu}{\gd^2}.
\end{equation}
A solution $u \in L^1_{loc}(\Gw)$ of the equation $-L_\mu u=0$ is called an \textit{$L_\mu$-harmonic function}.  Similarly, if
$$-L_\gm u\geq 0 \q \text{or}\q -L_\gm u\leq 0$$
we say that $u$ is \textit{$L_\gm$-superharmonic} or \textit{$L_\mu$-subharmonic}  respectively. If $\mu=0$ we shall just use the terms  harmonic, superharmonic, subharmonic.

Some problems involving the equation \eqref{N} with $\mu<1/4$ were studied by Bandle, Moroz and Reichel \cite{BMR}. They derived estimates of local $L_\mu$-subharmonic and superharmonic functions  and applied these results to study conditions for existence or nonexistence of large solutions of \eqref{N}. They also showed that the classical Keller -- Osserman estimate \cite{Ke, Oss} remains valid for \eqref{N}.


The condition $\gm<\frac{1}{4}$ is related to Hardy's inequality.


Denote by $C_H(\Gw)$ the best constant in Hardy's inequality, i.e.,
\bel{Hardy} C_H(\Gw)= \inf_{H^1_0(\Gw)} \frac{\int_\Gw |\nabla u|^2dx}{\int_\Gw (u/\gd)^2dx}.
\ee
By Marcus, Mizel and Pinchover \cite{MMP},  $C_H(\Gw) \in (0,\frac{1}{4}]$ and $C_H(\Gw)=\frac{1}{4}$ when $\Gw$ is convex. Furthermore the infimum is achieved if and only if $C_H(\Gw)<1/4$.
By Brezis and Marcus \cite{BrMa},
 for every $\mu< 1/4$ there exists a number $\gl_{\mu,1}$ such that
 $$\mu=\inf_{H^1_0(\Gw)}\frac{\int_\Gw (|\nabla u|^2-\gl_{\mu,1}u^2)dx}{\int_\Gw (u/\gd)^2dx}$$
 and the infimum is achieved. Thus $\gl_{\mu,1}$ is an eigenvalue of  $-L_\mu$; it is, in fact, a simple eigenvalue. We denote by $\vgf_{\mu,1}$ the corresponding positive eigenfunction normalized by
  $\int_\Gw(\vgf_{\mu,1}^2/\gd^2)dx=1$.

  If $\mu<C_H(\Gw)$ then $\gl_{\mu,1}$ is positive. Therefore, in this case, $\vgf_{\mu,1}$ is a positive supersolution of $-L_\mu$.
This fact and a classical result of Ancona \cite{An2} imply that the Martin boundary of $\Gw$ coincides with the Euclidean boundary and, if $K_\mu^\Gw$ denotes the Martin kernel for $-L_\mu$ in $\Gw$, the following theorem holds:
\medskip

 \noindent\textbf{Representation Theorem} \textit{For every $\nu\in \GTM^+(\bdw)$ the function
\begin{equation}\label{BBK}
    \BBK_\mu^\Gw[\nu](x):=\int_{\prt\Gw}K_\mu^{\Gw}(x,y)d\nu(y) \q\forall x\in \Gw
\end{equation}
is $L_\mu$-harmonic, i.e., $L_\mu \BBK_\mu^\Gw[\nu]=0$. Conversely, for every positive $L_\mu$-harmonic function $u$ there exists a unique measure $\nu\in \GTM^+(\prt\Gw)$ such that $u=\BBK_\mu^\Gw[\nu]$.}

\medskip

The measure $\nu$  such that $u=\BBK_\mu^\Gw[\nu]$ is called the \emph{$L_\mu$-boundary measure of $u$}. If $\mu=0$, $\nu$ is the measure boundary trace of $u$ (see \rdef{M-trace}). But if $\mu\in (0,C_H(\Gw))$ it can be shown that, for every $\nu\in \GTM^+(\bdw)$,  the measure boundary trace of $\BBK_\mu^\Gw[\nu]$ is zero.

\medskip

In the case $\mu=0$, the boundary value problem
\begin{equation}\label{PN0}\BAL
-\Gd u+|u|^{q-1}u&=0 \q \text{in }\Gw\\
u&=\nu \q \text{on } \prt\Gw
 \EAL
 \end{equation}
where $q>1$ and $\nu$ is either a finite measure or a positive  (possibly unbounded) measure,
has been studied by numerous authors. Following Brezis \cite{Br72}, if $\nu$ is a finite measure, a weak solution of \eqref{PN0} is defined as follows: $u$ is a solution of the problem if $u$ and $\gd |u|^q$ are integrable in $\Gw$ and
\bel{formN0} \int_{\Gw}{}(-u \Gd \zeta + |u|^{q-1}u \zeta)dx = -\int_{\prt \Gw}{}\frac{\prt \zeta}{\prt{\mathbf n}}d\gn \forevery \zeta \in C_0^2(\Gw) \ee
where ${\bf n}$ is the outer unit normal on $\prt\Gw$. Brezis proved that, if a solution exists then it is unique.
Gmira and V\'eron \cite{GV} showed that there exists  a critical exponent, $q_c:=\frac{N+1}{N-1}$, such that if $1<q<q_c$, \eqref{formN0} has a weak solution for every finite measure $\nu$ but, if $q\geq q_c$ there exists no positive solution with isolated point singularity.

Marcus and V\'eron \cite{MV1} proved that every positive solution of the equation
\begin{equation}\label{N0}
   -\Gd u+u^q=0
\end{equation}
possesses a boundary trace given by a positive measure $\nu$, not necessarily bounded. In the subcritical case the blow-up set of the trace is  a closed set. Furthermore they showed that, in this case, for every such positive measure $\nu$, the boundary value problem \eqref{PN0} has a unique solution.

 In the case $q=2$, $N=2$ this result was previously proved by Le Gall \cite{LG1} who introduced a probabilistic definition of the boundary trace.

 In the supercritical case the problem turned out to be much more challenging. It was studied by several authors using various techniques. The problem was studied by Le Gall, Dynkin, Kuznetsov, Mselati a.o. employing  mainly probabilistic methods. Consequently the results applied only to $1<q\leq 2$.  In parallel it was studied  by Marcus and Veron employing purely analytic methods that were not subject to the restriction $q\leq 2$. A complete classification of the positive solutions of \eqref{PN0} in terms of their behavior at the boundary was provided by Mselati \cite{Mse} for $q=2$, by Dynkin \cite{Dybook2} for $q_c\leq q\leq 2$ and finally by Marcus \cite{Ma} for every $q\geq q_c$. For details and related results we refer the reader to \cite{MVbook, MV4, MV3, AnMa, Dybook1} and the references therein.

In the case of equation \eqref{N} one is faced by the problem that, according to the classical definition of measure boundary trace, \textit{every positive $L_\mu$-harmonic function has measure boundary trace zero}. Therefore, in order to classify the positive solutions of \eqref{N} in terms of their behavior at the boundary, it is necessary to introduce a different notion of trace. As in the study of \eqref{N0}, we first consider the question of boundary trace for positive
$L_\mu$-harmonic or superharmonic functions.
\medskip

We recall the classical definition of measure boundary trace.

\bdef{M-trace} i) A sequence $\{D_n\}$ is a $C^2$ \emph{ exhaustion} of $\Gw$ if for every $n$, $D_n$ is of class $C^2$, $\ovl D_n \sbs D_{n+1}$ and $\cup_{n}D_n = \Gw$. If the domains are uniformly of class $C^2$ we say that $\{D_n\}$ is a uniform $C^2$ exhaustion.\medskip

\noindent ii) Let $u \in W^{1,p}_{loc}(\Gw)$ for some $p>1$. We say that $u$ possesses a \emph{measure boundary trace} on $\prt \Gw$  if there exists a finite measure $\gn$ on $\prt \Gw$ such that, for every uniform $C^2$ exhaustion $\{D_n\}$ and every $\vgf \in C(\ovl \Gw)$,
$$ \lim_{n \to \infty}\int_{\prt D_n} u|_{\prt D_n}\vgf dS = \myint{\prt \Gw}{}\vgf d\gn. $$
Here $u|_{D_n}$ denotes the Sobolev trace. The measure boundary trace of $u$ is denoted by $tr(u)$.
\es

\medskip

For $\gb>0$, denote
$$ \Gw_\gb=\{x \in \Gw: \gd(x)<\gb \},\; D_\gb=\{x \in \Gw: \gd(x)>\gb\}, \; \Gs_\gb=\{x \in \Gw: \gd(x)=\gb\}.$$

Put
\bel{alpha} \ga_\pm:=\myfrac{1}{2} \pm \sqrt{\myfrac{1}{4}-\gm}. \ee

It can be shown (see \rcor{2sides} below) that the classical \emph{measure boundary trace} of $\BBK_\mu^\Gw[\nu]$ is zero but there exist constants $C_1,C_2$ such that, for every $\nu\in \GTM(\prt\Gw)$,
\begin{equation}\label{2sides'}
  C_1 \norm{\gn}_{\GTM(\prt \Gw)} \leq    \rec{\gb^{\ga_-}} \int_{\Gs_\gb}\BBK_\mu^{\Gw}[\nu](x)dS_x \leq C_2\norm{\gn}_{\GTM(\prt \Gw)}
\end{equation}
for all $\gb\in (0,\gb_0)$ where $\gb_0>0$ depends only on $\Gw$. In view of this  we introduce the following definition of trace.
\bdef{nomtrace}
A positive function $u$ possesses a \emph{normalized boundary trace} if there exists a  measure $\gn\in\GTM^+(\prt\Gw)$ such that
\begin{equation}\label{tra}
   \lim_{\gb\to 0}\rec{\gb^{\ga_-}} \int_{\Gs_\gb}|u-\BBK_\mu^{\Gw}[\nu]|dS_x=0.
\end{equation}
The normalized boundary trace will be denoted by $tr^*(u)$.
\es

\Remark The notion of normalized boundary trace is well defined. Indeed, suppose that $\gn$ and $\gn'$ satisfy \eqref{tra}. Put $v=(\BBK_\gm^\Gw[\gn-\gn'])_+$ then $v$ is a nonnegative $L_\gm$-subharmonic function, $v\leq \BBK[\nu+\nu']$ and $tr^*(v)=0$. By \rprop{subhar}, $v=0$, i.e., $\BBK_\gm^\Gw[\gn-\gn'] \leq 0$. By interchanging the roles of $\gn$ and $\gn'$, we deduce that $\BBK_\gm^\Gw[\gn'-\gn] \leq 0$. Thus $\gn=\gn'$.
\medskip

Denote by $G_\mu^\Gw$ the Green function of $-L_\mu$ in $\Gw$ and, for every positive Radon measure $\tau$ in $\Gw$, put
$$\BBG_\mu^\Gw[\tau](x):=\int_\Gw G_\mu^\Gw(x,y)d\tau(y)$$

It can be shown that, if $\BBG_\mu^\Gw[\tau](x)<\infty$ for some $x\in \Gw$ then $\tau\in \GTM_{\gd^\gap}(\Gw)$. On the other hand, if $\tau\in \GTM^+_{\gd^\gap}(\Gw)$ then $\BBG_\mu^\Gw[\tau]$ is an $L_\mu$-harmonic function.
\medskip

 We begin with the study of the linear boundary value problem,

\begin{equation}\label{eqG}
\BAL
-L_\mu u&=\tau \q \text{in $\Gw$}\\
tr^*(u)&=\nu,
\EAL
\end{equation}
where $\nu\in \GTM^+(\prt\Gw)$ and $\tau\in \GTM^+_{\gd^\gap}(\Gw)$.
As usual we look for solutions $u\in L^1_{loc}(\Gw)$ and the equation is understood in the sense of distributions.
The representation theorem implies that if $\tau=0$ the problem has a unique solution, $u=\BBK_\mu^\Gw[\nu]$.

We list below our main results regarding this problem.
\medskip

\noindent\textbf{Proposition I. }

(i) \textit{If $u$ is a non-negative  $L_\mu$-harmonic function and $tr^*(u)=0$ then $u=0$.}
\smallskip

(ii) \textit{If $\tau\in \GTM^+_{\gd^\gap}(\Gw)$ then $\BBG_\mu^\Gw[\tau]$ has normalized trace zero. Thus $\BBG_\mu^\Gw[\tau]$ is a solution of \eqref{eqG}
with $\nu=0$.}
\smallskip

(iii) \textit{Let $u$ be a positive $L_\mu$-subharmonic function. If $u$ is dominated by an $L_\mu$-superharmonic function then $L_\mu u\in \GTM^+_{\gd^\gap}(\Gw)$ and $u$ has a normalized boundary trace. In this case $tr^*(u)=0$ if and only if $u\equiv0$.}
\smallskip

(iv) \textit{Let $u$ be a positive $L_\mu$-superharmonic function. Then there exist
$\nu\in \GTM^+(\bdw)$ and $\tau\in \GTM^+_{\gd^\gap}(\Gw)$ such that
\begin{equation}\label{rep1}
u=\BBG_\mu^\Gw[\tau]+ \BBK_\mu^\Gw[\nu].
\end{equation}
In particular, $u$ is an $L_\mu$-potential (i.e., $u$ does not dominate any positive $L_\mu$-harmonic function) if and only if $tr^*(u)=0$.
}

(v) \textit{For every $\nu\in \GTM^+(\bdw)$ and $\tau\in \GTM^+_{\gd^\gap}(\Gw)$, problem \eqref{eqG} has a unique solution. The solution is given by \eqref{rep1}.}
\medskip

Next we study the nonlinear boundary value problem,
\begin{equation}\label{PN}\BAL
-L_\mu u+u^q&=0 \q \text{in $\Gw$}\\
tr^*(u)&=\nu
\EAL
\end{equation}
where $\nu\in \GTM^+(\bdw)$.

\bdef{solN} (i) A positive solution of \eqref{N} is \emph{$L_\gm$-moderate} if it is dominated by an $L_\gm$-harmonic function.
\medskip

\noindent (ii) A positive function $u \in L^q_{loc}(\Gw)$ is a \emph{(weak) solution of \eqref{PN}} if  it satisfies the equation (in the sense of distributions) and has normalized boundary trace $\nu$.
\es

\bdef{XGw} Put
$$X(\Gw)=\{\zeta \in C^2(\Gw): \gd^{\ga_-}L_\gm\zeta \in L^\infty(\Gw), \; \gd^{-\ga_+}\zeta \in L^\infty(\Gw) \}. $$
A function $\gz\in X(\Gw)$ is called an \emph{admissible test function} for \eqref{PN}.
\es


Following are our main results concerning the nonlinear problem \eqref{PN}. Theorems A -- D apply to arbitrary exponent $q>1$.
\medskip

\noindent{\bf Theorem A.} {\it Assume that $0<\gm<C_H(\Gw)$, $q>1$. Let $u$ be a positive solution of \eqref{N}. Then the following statements are equivalent: \smallskip

\noindent{\bf i)} $u$ is $L_\gm$-moderate. \smallskip

\noindent{\bf ii)} $u$ admits a normalized boundary trace $\gn \in \GTM^+(\prt \Gw)$. In other words, $u$ is a solution of \eqref{PN}.\smallskip

\noindent{\bf iii)} $u \in L^q_{\gd^{\ga_+}}(\Gw)$ and
\bel{form} u + \BBG_\gm^\Gw[u^q] = \BBK_\gm^\Gw[\gn] \ee
where $\nu=tr^*(u)$. \smallskip

Furthermore, a positive function $u$  is a solution of \eqref{PN} if and only if
$u/\gd^{\ga_-} \in L^1(\Gw)$, $\gd^\gap u^q\in L^1(\Gw)$ and
\bel{intrep} \myint{\Gw}{} (- u L_\gm\zeta + u^q\zeta) dx = -\myint{\Gw}{}\BBK_\gm^\Gw[\gn]L_\gm \zeta dx \forevery \zeta \in X(\Gw). \ee
 }
\medskip

\noindent{\bf Theorem B.} {\it Assume $0<\gm<C_H(\Gw)$, $q>1$. \smallskip

\noindent {\sc i. Uniqueness.} For every $\gn \in \GTM^+(\prt \Gw)$, there exists at most one positive solution of \eqref{PN}. \smallskip

\noindent{\sc ii. Monotonicity.} Assume $\gn_i \in \GTM^+(\prt \Gw)$, $i=1,2$. Let $u_{\gn_i}$ be the unique solution of \eqref{PN} with $\gn$ replaced by $\gn_i$, $i=1,2$. If $\gn_1 \leq \gn_2$ then $u_{\gn_1} \leq u_{\gn_2}$.

\noindent{\sc iii. A-priori estimate.}  There exists a positive constant $c=c(N,\gm,\Gw)$ such that every  positive solution $u$ of \eqref{PN} satisfies,
\bel{estnormN} \norm{u}_{L^1_{\gd^{-\ga_-}}(\Gw)} + \norm{u}_{L^q_{\gd^{\ga_+}}(\Gw)} \leq c\norm{\gn}_{\GTM(\prt \Gw)}.
\ee
} \medskip

\noindent{\bf Theorem C.}  {\it Assume $0<\gm<C_H(\Gw)$, $q>1$. If $\nu\in \GTM^+(\bdw)$ and $\BBK_\mu^\Gw[\nu]\in L^q_{\gd^\gap}(\Gw)$ then there exists a unique solution of the boundary value problem \eqref{PN}.}
\medskip

\noindent{\bf Corollary C1.} \textit{For every positive function $f \in L^1(\prt \Gw)$  \eqref{PN} with $\nu=f$ admits a unique positive solution.} \medskip
\medskip

\noindent{\bf Theorem D.} {\it Assume $0<\gm<C_H(\Gw)$, $q>1$.  If $u$ is a positive solution of \eqref{PN} with $\gn\in \GTM^+(\bdw)$ then
\bel{behavior} \lim_{x \to y}\myfrac{u(x)}{\BBK_\gm^\Gw[\gn](x)}=1 \q \text{non-tangentially},\; \gn\text{-a.e. on } \prt \Gw. \ee
} \medskip

Let
\begin{equation}\label{qmc}
 q_{\gm,c}:=\frac{N+\ga_+}{N-1-\ga_-}.
\end{equation}
In the next two results we show, among other things, that $q_{\gm,c}$ is the \emph{critical exponent} for \eqref{PN}. This means that, if $1<q<q_{\mu,c}$ then problem \eqref{PN} has a unique solution for every measure $\nu\in \GTM^+(\bdw)$ but, if $q\geq q_{\mu,c}$ then the problem has no solution for some measures $\nu$, e.g. Dirac measure.

In Theorem E  we consider the subcritical case $1<q<q_{\gm,c}$ and in Theorem F the supercritical case.
\medskip

\noindent{\bf Theorem E.} \textit{ Assume $0<\gm<C_H(\Gw)$ and $1<q<q_{\gm,c}$. Then:}

\noindent{\sc i. Existence and uniqueness.}\hskip 2mm\textit{For every $\gn \in \GTM^+(\prt \Gw)$ \eqref{PN} admits a unique positive solution $u_\gn$.}\medskip

\noindent{\sc ii. Stability.} \hskip 2mm
\textit{If $\{\gn_n\}$ is a sequence of measures in $\GTM^+(\prt \Gw)$  weakly convergent to $\gn \in \GTM^+(\prt \Gw)$ then  $u_{\gn_n}\to u_\gn$  in $L^1_{\gd^{-\ga_-}}(\Gw)$ and in $L^q_{\gd^{\ga_+}}(\Gw)$.}
 \medskip

\noindent{\sc iii. Local behavior.}\hskip 2mm \textit{Let $\gn=k\gd_y$, where $k>0$ and $\gd_y$ is the Dirac measure concentrated at $y \in \prt \Gw$. Then, under the assumptions of Theorem E, the unique solution of \eqref{PN}, denoted by $u_{k\gd_y}$, satisfies
\bel{dirac} \lim_{x \to y}\myfrac{u_{k\gd_y}(x)}{K_\gm^\Gw(x,y)}=k. \ee
} \medskip

\noindent{\it Remark.} \hskip 2mm Note that in part (iii) we have `uniform convergence' not just `non-tangential convergence' as in Theorem D. \medskip

\noindent{\bf Theorem F.} {\it Assume $0<\gm<C_H(\Gw)$ and $q \geq q_{\gm,c}$. Then for every $k>0$ and $y \in \prt \Gw$, there is no positive solution of \eqref{N} with normalized boundary trace $k \gd_y$.} \medskip

 In the first part of the paper we  study properties of positive $L_\mu$-harmonic functions and the boundary value problem \eqref{eqG}.  In the second part, these results are applied to a study of the corresponding boundary value problem for the nonlinear equation \eqref{N}. These results yield a complete classification of the positive moderate solutions of \eqref{N} in the subcritical case. They also provide a framework for the study of positive solutions of \eqref{N} that may blow up at some parts of the boundary. The existence of such solutions in the subcritical case has been studied (by different methods) in \cite{BMM}. Corresponding boundary value problems - including a study of solutions with strong isolated singularities - will be presented in a forthcoming paper
 \cite{Msing}.

The main ingredients used in this paper are: the Representation Theorem previously stated and other basic results of potential theory (see \cite{An1}), a sharp estimate of the Green kernel of $-L_\mu$ due to Filippas, Moschini and Tertikas \cite{FMT}, estimates for convolutions in weak $L^p$ spaces (see \cite[Section 2.3.2]{MVbook}) and the comparison principle obtained in \cite{BMR}.

\section{The linear equation}
Throughout this paper we assume that $0<\gm<C_H(\Gw)$.

\subsection{Some potential theoretic results}
We denote by $\GTM_{\gd^\ga}(\Gw)$, $\ga \in \BBR$, the space of Radon
measures $\gt$ on $\Gw$ satisfying $\int_{\Gw}\gd^\ga(x) d|\gt|<\infty$ and by
$\GTM_{\gd^\ga}^+(\Gw)$ the positive cone of  $\GTM_{\gd^\ga}(\Gw)$.  When
$\ga=0$, we use the notation  $\GTM(\Gw)$ and $\GTM^+(\Gw)$. We also denote by
$\GTM(\prt \Gw)$ the space of finite Radon measures on $\prt \Gw$ and by
$\GTM^+(\prt \Gw)$ the positive cone of $\GTM(\prt \Gw)$.

Let $D$ be a $C^2$ domain such that $D \sbs \sbs \Gw$ and $h \in L^1(\prt\Gw)$. Denote by $\BBS_\gm(D,h)$ the solution of the problem
\bel{D} \left\{ \BA{lll} -L_\gm u &= 0 \qq &\text{in } D \\
\phantom{--,}
u &= h &\text{on } \prt D.
\EA \right. \ee

\blemma{sup1} Let $u$ be $L_\gm$-superharmonic in $\Gw$ and $D$ be a $C^2$
domain such that $D \sbs \sbs \Gw$. Then $u \geq \BBS_\gm(D,u)$ a.e. in $D$.
\es
\proof Since $u$  is $L_\gm$-superharmonic in $\Gw$, there exists $\gt \in
\GTM^+(\Gw)$ such that $-L_\gm u =\gt$. Let $v$ be the solution of
\bel{D1} \left\{ \BA{lll} -L_\gm v &= \gt \qq &\text{in } D \\
\phantom{--,}
v &= 0 &\text{on } \prt D
\EA \right. \ee
and $w=\BBS_\gm(D,u)$. Then $w \geq 0$ and $u|_D=v+w \geq v$. \qed

\blemma{largest} Let $u$ be a nonnegative $L_\gm$-superharmonic and $\{D_n\}$
be a $C^2$ exhaustion of $\Gw$. Then
$$ \hat u:=\lim_{n \to \infty}\BBS_\gm(D_n,u) $$
exists and is the largest $L_\gm$-harmonic function dominated by $u$.
\es
\proof By \rlemma{sup1}, $\BBS_\gm(D_n,u) \leq u|_{D_n}$, hence the sequence
$\{ \BBS_\gm(D_n,u)\}$ is decreasing. Consequently, $\hat u$ exists and is an
$L_\gm$-harmonic function dominated by $u$. Next, if $v$ is an
$L_\gm$-harmonic function dominated by $u$ then $v \leq \BBS_\gm(D_n,u)$ for
every $n \in \BBN$. Letting $n \to \infty$ yields $v \leq \hat u$. \qed
\bdef{potential} A nonnegative $L_\gm$-superharmonic function is called an
{\it $L_\gm$-potential} if its largest $L_\gm$-harmonic minorant is zero. \es
As a consequence of \rlemma{largest}, we obtain
\blemma{potential} Let $u_p$ be a nonnegative $L_\gm$-superharmonic function
in $\Gw$. If for some $C^2$ exhaustion $\{D_n \}$ of $\Gw$,
\bel{lim0} \lim_{n \to \infty}\BBS_\gm(D_n,u_p) = 0, \ee
then $u_p$ is an $L_\gm$-potential in $\Gw$. Conversely, if $u_p$ is an
$L_\gm$-potential, then \eqref{lim0} holds for every $C^2$ exhaustion
$\{D_n\}$ of $\Gw$.
\es
\medskip
For easy reference we quote below the Riesz decomposition theorem (see \cite{An1}).

\bth{Ri} Every nonnegative $L_\gm$-superharmonic function $u$ in $\Gw$ can be
written in a unique way in the form $u=u_p+u_h$ where $u_p$ is an
$L_\gm$-potential and $u_h$ is a nonnegative $L_\gm$-harmonic function in
$\Gw$.
\es

The next result is a consequence of the Fatou convergence theorem  \cite[Theorem 1.8]{An1} and the following well-known fact:  if a function satisfies the Harnack inequality, fine convergence at the boundary (in the sense of \cite{An1}) implies non-tangential convergence.
\bth{pot} Let $u_p$ be a positive $L_\gm$-potential and $u$ be a positive
$L_\gm$-harmonic function. Assume that $\frac{u_p}{u}$ satisfies the Harnack
inequality. Then
$$ \lim_{x \to y}\myfrac{u_p(x)}{u(x)}=0 \q \text{non-tangentially, }
\gn\text{-a.e. on } \prt \Gw $$
where $\gn$ is the $L_\gm$-boundary measure of $u$.
\es

\subsection{The action of the Green and Martin kernels on spaces of measures}

From \cite{An2}, for every $y \in \prt \Gw$, there exists a positive
$L_\gm$-harmonic function in $\Gw$ which vanishes on $\prt \Gw \sms \{y\}$.
When normalized, this function is unique. We choose a fixed reference point
$x_0$ in $\Gw$ and denote by $K_{\gm,y}^{\Gw}$ this $L_\gm$-harmonic function,
normalized by $K_{\gm,y}^{\Gw}(x_0)=1$. The function
$K_\gm^{\Gw}(\cdot,y)=K_{\gm,y}^{\Gw}(\cdot)$ is the $L_\gm$-Martin kernel in
$\Gw$, normalized at $x_0$.

For $\gn\in \GTM(\bdw)$ denote
$$ \BBK_\gm^\Gw[\gn](x)=\myint{\prt \Gw}{}K_\gm^\Gw(x,y)d\gn(y).$$

Let $G_\gm^\Gw$ be the Green kernel for the operator $L_\gm$ in $\Gw \times
\Gw$.
If $\gt \in \GTM_{\gd^{\ga_+}}(\Gw)$ then
$$\BBG_\mu^\Gw[\tau](x):=\myint{\Gw}{}G_\gm^\Gw(x,y)d\gt(y)<\infty \q\text{a.e. in } \Gw. $$

Denote by $G^\Gw:=G_0^\Gw$ and $P^\Gw:=P_0^\Gw$ the Green and Poisson kernels
respectively of $-\Gd$ in $\Gw$. We recall that (see, e.g., \cite{MVbook})
$$ G^\Gw(x,y) \sim \min\{\gd(x),\gd(y)\}|x-y|^{1-N} \forevery x,y \in \Gw, x
\ne y, $$
$$ P^\Gw(x,y) \sim \gd(x)|x-y|^{-N} \forevery x \in \Gw, y \in \prt \Gw, $$
where the notation $f \sim g$ means there exists a postive constant $c$ such
that $c^{-1}f < g < cf$. By \cite[Theorem 4.11]{FMT}, for every $x,y \in \Gw,
x \ne y$,
\bel{Gmu} G^\Gw_\gm(x,y) \sim
\min\left\{\abs{x-y}^{2-N},\gd(x)^{\ga_+}\gd(y)^{\ga_+}\abs{x-y}^{2\ga_-
-N}\right\} \ee
Since
$$ K_\gm^\Gw(x,y):=\lim_{z \to y}\myfrac{G^\Gw_\gm(x,z)}{G^\Gw_\gm(x_0,z)}
\forevery x \in \Gw $$
it follows from \eqref{Gmu} that
\bel{Pmu} K_\gm^\Gw(x,y) \sim \gd(x)^{\ga_+}|x-y|^{2\ga_- -N} \forevery x \in
\Gw, y \in \prt \Gw. \ee
As a consequence,
\bel{equi} \myfrac{K_\gm^\Gw(x,y)}{\gd(x)^{\ga_-}} \sim
\myfrac{\gd(x)}{|x-y|^N}\left( \myfrac{|x-y|}{\gd(x)} \right)^{2\ga_-} \sim
P^\Gw(x,y) \left( \myfrac{|x-y|}{\gd(x)} \right)^{2 \ga_-}. \ee

\Remark (i) By \eqref{Gmu}, $G_\mu^\Gw(x_0,x)\sim \gd(x)^{\ga_+}$ as $\gd(x)\to 0$. It is well known that the first eigenfunction behaves, near the boundary, in the same way as the Green function. Therefore,
\begin{equation}\label{vgf1}
   c^{-1}\gd(x)^{\ga_+}\leq \vgf_{\mu,1}\leq c \gd(x)^{\ga_+}.
\end{equation}
(ii) We note that  $P^{\BBR^N_+}(x,0)=x_1|x|^{-N}$ and
$P_\gm^{\BBR^N_+}(x,0)=x_1^{\ga_+}|x|^{2\ga_+-N}$. \medskip
\medskip

Denote  $L^p_w(\Gw;\tau)$, $1 \leq p < \infty$, $\tau \in \GTM^+(\Gw)$, the
weak $L^p$ space defined as follows: a measureable function $f$ in $\Gw$
belongs to this space if there exists a constant $c$ such that
\bel{distri} \gl_f(a;\tau):=\tau(\{x \in \Gw: |f(x)|>a\}) \leq ca^{-p},
\forevery a>0. \ee
The function $\gl_f$ is called the distribution function of $f$ (relative to
$\tau$). For $p \geq 1$, denote
$$ L^p_w(\Gw;\tau)=\{ f \text{ Borel measurable}:
\sup_{a>0}a^p\gl_f(a;\tau)<\infty\} $$
and
\bel{semi}
\norm{f}^*_{L^p_w(\Gw;\tau)}=(\sup_{a>0}a^p\gl_f(a;\tau))^{\frac{1}{p}}. \ee
The $\norm{.}_{L^p_w(\Gw;\tau)}$ is not a norm, but for $p>1$, it is
equivalent to the norm
\bel{normLw} \norm{f}_{L^p_w(\Gw;\tau)}=\sup\left\{
\frac{\int_{\gw}|f|d\tau}{\tau(\gw)^{1/p'}}:\gw \sbs \Gw, \gw \text{
measurable }, 0<\tau(\gw)<\infty \right\}. \ee
More precisely,
\bel{equinorm} \norm{f}^*_{L^p_w(\Gw;\tau)} \leq \norm{f}_{L^p_w(\Gw;\tau)}
\leq \myfrac{p}{p-1}\norm{f}^*_{L^p_w(\Gw;\tau)}. \ee
Notice that, for every $\ga>-1$,
$$L_w^p(\Gw;\gd^{\ga}dx) \sbs L_{\gd^{\ga}}^{r}(\Gw) \forevery r \in [1,p).
$$

For every $x \in \prt \Gw$, denote by ${\bf n}_x$ the outward unit normal
vector to $\prt \Gw$ at $x$.

The following is a well-known geometric property of $C^2$ domains.
\bprop{beta0} There exists $\gb_0>0$ such that

\noindent i) For every point $x \in \ovl \Gw_{\gb_0}$, there exists a unique
point $\gs_x \in \prt \Gw$ such that $|x -\gs_x|=\gd(x)$. This implies
$x=\gs_x - \gd(x){\bf n}_{\gs_x}$.

\noindent ii) The mappings $x \mapsto \gd(x)$ and $x \mapsto \gs_x$ belong to
$C^2(\ovl \Gw_{\gb_0})$ and $C^1(\ovl \Gw_{\gb_0})$ respectively. Furthermore,
$\lim_{x \to \gs(x)}\nabla \gd(x) = - {\bf n}_x$.

\es

By combining \eqref{Gmu}, \eqref{Pmu} and the estimates of \cite[Lemma
2.3.2]{MVbook}, we obtain
\bprop{GP} There exist constants $c_i$ depending only on $N,\gm,\gb, \Gw$ such that,
\bel{estG} \norm{\BBG^\Gw_\gm[\gt]}_{L_w^{\frac{N+\gb}{N-2}}(\Gw,\gd^\gb)}
\leq c_1\norm{\gt}_{\GTM(\Gw)},\q \forall \gt \in \GTM(\Gw),\;\gb>-1, \ee
\bel{estG1}
\norm{\BBG^\Gw_\gm[\gt]}_{L_w^{\frac{N+\gb}{N-2\ga_-}}(\Gw,\gd^{\gb-\ga_+})} \leq
c_1\norm{\gt}_{\GTM_{\gd^{\ga_+}}(\Gw)}, \q\forall \gt \in
\GTM_{\gd^{\ga_+}}(\Gw),\;\gb>-2\gan, \ee
\bel{estP}
\norm{\BBK_\gm^\Gw[\gn]}_{L_w^{\frac{N+\gb}{N-1-\ga_-}}(\Gw,\gd^\gb)} \leq
c_2\norm{\gn}_{\GTM(\prt \Gw)},\q \forall \gn \in \GTM(\prt \Gw),\;  \gb>-1. \ee
\es
\proof We assume that $\gt$ is positive; otherwise we replace $\gt$ by
$|\gt|$. We consider $\gt$ as a positive measure in $\BBR^N$ by extending
$\gt$ by zero outside of $\Gw$. For $\ga \in (0,N)$, denote
$\Gg_\ga(x)=|x|^{\ga-N}$. By \cite[Lemma 2.3.3]{MVbook} (see inequality (2.3.17)),
\begin{equation}\label{convolution}
   \norm{\Gg_\ga*\gt}_{L_w^{\frac{N+\gb}{N-\ga}}(\Gw,\gd^{\gb})} \leq
c\norm{\gt}_{\GTM(\Gw)} \forevery \gb>\max\{-1,-\ga\}
\end{equation}
where $c=c(N,\ga,\gb, diam(\Gw))$. By $(\ref{Gmu})$,
$$G^\Gw_\gm(x,y) \leq c\min\{\Gg_2(x-y), \gd(x)^{\ga_+}\gd(y)^{\ga_+}\Gg_{2\ga_-}(x-y)\}.$$
Hence, by \eqref{convolution},
$$ \BAL \norm{\BBG^\Gw_\gm[\gt]}_{L_w^{\frac{N+\gb}{N-2}}(\Gw,\gd^{\gb})}
&\leq c\norm{\Gg_2 * \gt}_{L_w^{\frac{N+\gb}{N-2}}(\Gw,\gd^\gb)}\\
&\leq c'\norm{\gt}_{\GTM(\Gw)} \forevery\gb>-1,\\[2mm]
 \norm{\BBG^\Gw_\gm[\gt]}_{L_w^{\frac{N+\gb}{N-2\ga_-}}(\Gw,\gd^{\gb-\gap})} &\leq
c\norm{\Gg_{2\ga_-} * (\gd^{\ga_+}\gt)}_{L_w^{\frac{N+\gb}{N-2\ga_-}}(\Gw,\gd^{\gb})}\\
 & \leq c\norm{\gt}_{\GTM_{\gd^{\ga_+}}(\Gw)} \forevery \gb>-2\ga_-.
\EAL $$
Next we extend  $\gn$ by zero outside $\prt \Gw$ and observe that,
by $(\ref{Pmu})$, $K_\gm^\Gw(x,y) \leq c\Gg_{1+\ga_-}(x-y)$. Hence
$\BBK_\gm^\Gw[\gn] \leq c\Gg_{1+\ga_-}*\gn$ and by \eqref{convolution},
$$ \norm{\BBK_\gm^\Gw[\gn]}_{L_w^{\frac{N+\gb}{N-1-\ga_-}}(\Gw,\gd^\gb)} \leq
c\norm{\Gg_{1+\ga_-} * \gn}_{L_w^{\frac{N+\gb}{N-1-\ga_-}}(\Gw,\gd^\gb)} \leq
c\norm{\gn}_{\GTM(\prt \Gw)} \forevery \gb>-1. $$ \qed
\medskip

\bcor{convergence}  Let $\gb>-1$.

\noindent(i)\hskip 2mm  If $\{\gn_n\} \sbs \GTM^+(\prt \Gw)$
converges weakly to $\gn \in \GTM^+(\prt \Gw)$ then $\{\BBK_\gm^\Gw[\gn_n]\}$
converges to $\BBK_\gm^\Gw[\gn]$ in $L^p_{\gd^\gb}(\Gw)$ for every $p$ such that  $1 \leq p <\frac{N+\gb}{N-1-\ga_-}$.

\noindent(ii)\hskip 2mm If $\{\tau_n\} \sbs \GTM^+(\Gw)$
converges weakly (relative to $C_0(\bar\Gw)$) to $\tau \in \GTM^+(\Gw)$ then $\{\BBG_\gm^\Gw[\tau_n]\}$
converges to $\BBG_\gm^\Gw[\tau]$ in $L^p_{\gd^\gb}(\Gw)$ for every $p$ such that  $1 \leq p <\frac{N+\gb}{N-2}$.
\es
\proof We prove the first statement. The second is proved in the same way.

Since $K^\Gw_\gm(x,.) \in C(\prt \Gw)$ for every $x \in \Gw$,
$\{\BBK_\gm^\Gw[\gn_n]\}$ converges to $\BBK_\gm^\Gw[\gn]$ every where in
$\Gw$. By Holder inequality and $(\ref{estP})$, we deduce that
$\{(\BBK_\gm^\Gw[\gn_n])^p\}$ is equi-integrable w.r.t. $\gd^\gb dx$ for any
$1\leq p < \frac{N+\gb}{N-1-\ga_-}$. By Vitali's theorem, $\BBK_\gm^\Gw[\gn_n]
\to \BBK_\gm^\Gw[\gn]$ in $L^p_{\gd^\gb}(\Gw)$. \qed

\subsection{Estimates related to the normalized trace}
\bprop{extrace} There exist positive constants $C_1,C_2$ such that, for every
$\gb\in (0,\gb_0)$,
\begin{equation}\label{tr*P}
  C_1\gb^{\ga_-}\leq \int_{\Gs_\gb}K_\gm^\Gw(x,y)dS_x \leq C_2 \gb^{\ga_-}
  \q\forall y\in \prt\Gw.
\end{equation}
The constants $C_1,C_2$ depend on $N,\Gw,\gm$ but not on $y$.

Furthermore, for every $r_0>0$,
\begin{equation}\label{tr*P0}
\lim_{\gb\to0} \frac{1}{\gb^{\ga_-}}\int_{\Gs_\gb\sms
B_{r_0}(y)}K_\gm^\Gw(x,y)dS_x=0 \q\forall y\in \prt\Gw.
\end{equation}
For $r_0$ fixed, the rate of convergence is independent of $y$.
\es

\proof By \eqref{Pmu},
\begin{equation}\label{estP0}
 \frac{1}{\gb^{\ga_-}}\int_{\Gs_\gb \sms B_{r_0}(y)}K_\gm^\Gw(x,y)dS_x\leq
 c\gb^{\ga_+-\ga_-}.
\end{equation}
This implies \eqref{tr*P0}.

For the next estimate it is convenient to assume that the coordinates are
placed so that $y=0$ and the tangent hyperplane to $\prt\Gw$ at $0$ is $x_N=0$
with the $x_N$ axis pointing into the domain. For $x\in \BBR^N$ put
$x'=(x_1,\cdots,x_{N-1})$. Pick $r_0\in (0,\gb_0)$  sufficiently small
(depending only on the $C^2$ characteristic of $\Gw$) so that
$$\frac{1}{2}(|x'|^2+\gd(x)^2)\leq |x|^2\q \forall x\in \Gw\cap B_{r_0}(0).$$
Then, if $x\in \Gs_{\gb}\cap B_{r_0}(0)=:\Gs_{\gb,0}$,
$$\frac{1}{4}( |x'| +\gb)\leq |x|.$$
This inequality and \eqref{Pmu} imply,
\begin{equation*}\BAL
 \int_{\Gs_{\gb,0}} K_\mu^\Gw(x,0)dS_x&\leq c_0\gb^{\ga_+}\int_{\Gs_{\gb,0}} (|x'|
 +\gb)^{2\ga_--N}dS_x\\
 &\leq c_1\gb^{\ga_+}\int_{|x'|<r_0} (|x'| +\gb)^{2\ga_--N}dx'\\
 &\leq c_2\gb^{\ga_+}\int_0^{r_0}(t+\gb)^{2\ga_--2}dt\\
 &<c_2\gb^{\ga_-}\int_1^{\infty}
 \tau^{-2\ga_+}d\tau=\frac{c_2}{2\ga_+-1}\gb^{\ga_-}.
 \EAL\end{equation*}
Thus, for $\gb<r_0$,
\begin{equation}\label{Pgb2}
 \frac{1}{\gb^{\ga_-}}\int_{\Gs_{\gb,0}} K_\mu^\Gw(x,0)dS_x \leq
 \frac{c_2}{2\ga_+-1}.
\end{equation}

Estimates \eqref{estP0} and \eqref{Pgb2} imply the second estimate in
\eqref{tr*P}. The first estimate in \eqref{tr*P} follows from \eqref{equi}. \qed

Since \eqref{tr*P} holds uniformly w.r. to $y\in \prt\Gw$,  an application of
Fubini's yields the following.

\bcor{2sides} For every $\gn \in \GTM^+(\prt \Gw)$,
\begin{equation}\label{2sides}
  \BAL C_1 \norm{\gn}_{\GTM(\prt \Gw)} &\leq \liminf_{\gb \to 0}
  \myint{\Gs_\gb}{}\myfrac{\BBK_\gm^\Gw[\gn]}{\gd(x)^{\ga_-}}dS_x \\
&\leq  \limsup_{\gb \to 0}
\myint{\Gs_\gb}{}\myfrac{\BBK_\gm^\Gw[\gn]}{\gd(x)^{\ga_-}}dS_x\leq
C_2\norm{\gn}_{\GTM(\prt \Gw)} \EAL
\end{equation}
with $C_1, C_2$ as in \eqref{tr*P}.
\es

\bprop{Gmt} If $\tau\in \GTM_{\gd^{\ga_+}}(\Gw)$
then
\begin{equation}\label{Gmt}
 tr^*(\BBG_\mu^\Gw[\tau])=0
\end{equation}
and, for $0<\gb<\gb_0$,
\begin{equation}\label{Gmt'}
   \rec{\gb^\gan} \int_{\Gs_\gb}\BBG_\mu^\Gw[\tau]dS_x\leq c\int_\Gw\gd^{\ga_+}d|\tau|,
\end{equation}
where $c$ is a constant depending on $\mu,\Gw$.
\es

\proof We may assume that $\tau>0$. Denote $v:=\BBG_\mu^\Gw[\tau].$ We start with the proof of \eqref{Gmt'}.

By Fubini's theorem and \eqref{Gmu},
\[\BAL
\int_{\Gs_\gb}vdS_x\leq c\Big(&\int_\Gw\int_{\Gs_\gb\cap B_{\frac{\gb}{2}}(y)}|x-y|^{2 -N} dS_x\,d\tau(y)\\
+\gb^{\gap}&\int_\Gw\int_{\Gs_\gb\sms B_{\frac{\gb}{2}}(y)}|x-y|^{2\gan -N} dS_x\,\gd^{\gap}(y)d\tau(y)\Big)=I_1(\gb)+I_2(\gb).
\EAL\]
Note that, if $x\in\Gs_\gb$ and $|x-y|\leq \gb/2$ then $\gb/2\leq \gd(y)\leq 3\gb/2$. Therefore
$$\BAL
I_1(\gb)&\leq c_1\int_{\Gs_\gb\cap B_{\frac{\gb}{2}}(y)}|x-y|^{2-\gap -N}dS_x\int_\Gw \gd(y)^\gap\, d\tau(y)\\
&\leq c'_1\int_0^{\gb/2}r^{2-\gap -N}r^{N-2}dr\,\int_\Gw \gd(y)^\gap\, d\tau(y)\\
&\leq c''_1\gb^\gan\,\int_\Gw \gd(y)^\gap\, d\tau(y)
\EAL$$
and
$$I_2(\gb)\leq c_2\gb^\gap\int_{\gb/2}^\infty r^{2\gan -N}r^{N-2}dr\int_\Gw \gd(y)^\gap\,d\tau=c'_2\gb^\gan \int_\Gw \gd(y)^\gap\,d\tau.$$
This implies \eqref{Gmt'}.

Given $\ge\in (0,\norm{\tau}_{\GTM_{\gd^\gap}(\Gw)})$ and $\gb_1\in (0,\gb_0)$ put $\tau_1=\tau\chi_{_{\bar D_{\gb_1}}}$ and $\tau_2=\tau-\tau_1$. Pick $\gb_1=\gb_1(\ge)$ such that
\begin{equation}\label{tau2}
 \int_{\Gw_{\gb_1}}\gd(y)^\gap\,d\tau\leq\ge.
\end{equation}
 Thus the choice of $\gb_1$ depends on the rate at which $\int_{\Gw_\gb}\gd^\gap\,d\tau $ tends to zero as $\gb\to 0$.

Put $v_i=\BBG_\mu^\Gw[\tau_i]$.
Then, for $0<\gb<\gb_1/2$,
$$\int_{\Gs_\gb}v_1\,dS_x\leq c_3\gb^{\gap}\gb_1^{2\gan-N}\int_\Gw \gd^{\gap}(y)d\tau_1(y).$$
Thus,
\begin{equation}\label{Gmt1}
   \lim_{\gb\to 0}\rec{\gb^\gan} \int_{\Gs_\gb}v_1\,dS_x=0.
\end{equation}
On the other hand, by \eqref{Gmt'} and \eqref{tau2},
\begin{equation}\label{Gmt2}
   \rec{\gb^\gan} \int_{\Gs_\gb}v_2\,dS_x\leq c\ge \forevery \gb<\gb_0.
\end{equation}
This implies that $tr^*(v)=0$.
\qed
\medskip

 \Remark
 It is well-known that  $u$ is an $L_\mu$-potential if and only if there exists a positive measure $\tau$ in $\Gw$ such that $u=\BBG_\mu^\Gw[\tau]$ (see e.g.\cite[Theorem 12]{An1}). The estimate  \eqref{Gmu} implies that if $\BBG_\mu^\Gw[\tau]\not\equiv\infty$ then $\tau\in \GTM_{\gd^{\gap}}(\Gw)$. Therefore as a consequence of the previous proposition:

 \bcor{potential}A positive $L_\mu$-superharmonic function $u$ is a potential if and only if $tr^*(u)=0$.
 \es

\bprop{subhar} Let $w$ be a non-negative $L_\mu$-subharmonic function.
 If $w$ is dominated by an $L_\mu$-superharmonic function then $L_\mu w\in \GTM^+_{\gd^\gap}(\Gw)$ and  $w$ has a normalized boundary trace $\nu\in \GTM^+(\bdw)$. If, in addition, $tr^*(w)=0$ then $w=0$.
\es
\proof  The first assumption implies that there exists a positive Radon measure $\gl$ in $\Gw$ such that $-L_\mu w=-\gl$.

First assume that $\gl\in \GTM_{\gd^\gap}(\Gw)$. Then $v:=w+\BBG_\mu^\Gw[\gl]$ is a non-negative $L_\mu$-harmonic function and consequently, by the representation theorem, $v=\BBK_\mu^\Gw[\nu]$ for some $\nu\in \GTM^+(\bdw)$. By \rprop{Gmt}, $tr^*(w)=\nu$. If $\nu=0$ then $v=0$ and therefore $w=0$. Now let us drop the assumption on $\gl$.

Let $v_\gb$ be the unique solution of the boundary value problem,
$$-L_\mu v_\gb=-\gl_\gb \;\txt{in}\; D_\gb,\q v_\gb=h_\gb \;\txt{on}\;\prt D_\gb$$
where $\gl_\gb$ is the restriction of $\gl$ to $D_\gb$ and $h_\gb$ is the restriction of $w$ to $\prt D_\gb$. (The uniqueness follows from \cite[Lemma 2.3]{BMR}.) The uniqueness implies that $v_\gb=w\lfloor_{D_\gb}$. By assumption there exists a positive $L_\mu$-superharmonic function, say $V$, such that $w\leq V$. Hence
$$w+\BBG_\mu^{D_\gb}[\gl_\gb]= P_\mu^{D_\gb}[h_\gb]\leq P_\mu^{D_\gb}[V\lfloor_{\prt D_\gb}]\leq V.$$
This implies that  $\BBG_\mu^\Gw[\gl]=\lim_{\gb\to 0}\BBG_\mu^{D_\gb}[\gl_\gb]<\infty$. For fixed $x\in \Gw$, $G_\mu^\Gw(x,y)\sim\gd(y)^\gap$. Therefore the finiteness of $\BBG_\mu^{\Gw}[\gl]$ implies that $\gl\in \GTM_{\gd^\gap}(\Gw)$. By the first part of the proof $w$ has a normalized trace.\qed
\medskip

\noindent\textit{Remark.} See \rprop{subhar'} below for a complementary result.

\subsection{Test functions}\label{S:test} Denote
$$X(\Gw)=\{\zeta \in C^2(\Gw): \gd^{\ga_-}L_\gm\zeta \in L^\infty(\Gw),
\gd^{-\ga_+}\zeta \in L^\infty(\Gw) \}. $$

\bprop{testgrad} For any $\zeta \in X(\Gw)$, $\gd^{\ga_-}|\nabla \zeta| \in
L^\infty(\Gw)$.
\es
\proof Let $\zeta \in X(\Gw)$ then there exist a positive constant $c_1$ and a
function $f \in L^\infty(\Gw)$ such that $|\zeta| \leq c_1 \gd^{\ga_+}$ and
$$ -L_\gm \zeta = \gd^{-\ga_-}f. $$
Take arbitrary point $x_* \in \Gw_{\gb_0}$ and put $d_*=\frac{1}{2}\gd(x_*)$,
$y_*=\frac{1}{d_*}x_*$, $\zeta_*(y)=\zeta(d_*y)$
for $y \in \frac{1}{d_*}\Gw_{d_*}$. Note that if $x \in B_{d_*}(x_*)$ then
$y=\frac{1}{d_*}x \in B_1(y_*)$ and $1 \leq \dist(y,\prt (\frac{1}{d_*}\Gw_{d_*})) \leq 3$. In
$B_1(y_*)$,
$$ - \Gd \zeta_* - \myfrac{\gm}{\dist(y,\prt (\frac{1}{d_*}\Gw_{d_*}))^2}\zeta_* =
d_*^{2-\ga_-}\dist(y,\prt (\frac{1}{d_*}\Gw_{d_*}))^{-\ga_-}f(d_*y). $$
By local estimate for elliptic equations \cite[Theorem 8.32]{GT}, there exists
a positive constant $c_2=c_2(N,\gm)$ such that
$$ \max_{B_{\frac{1}{2}}(y_*)}\abs{\nabla \zeta_*} \leq
c_2[\max_{B_1(y_*)}\abs{\zeta_*}+\max_{B_1(y_*)}(d_*^{2-\ga_-}\dist(y,\prt (\frac{1}{d_*}\Gw_{d_*}))^{-\ga_-}\abs{f(d_*y)}]. $$
This implies
$$d_*\abs{\nabla \zeta (x_*)} \leq c_3(\gd(x_*)^{\ga_+} +
\norm{f}_{L^\infty(\Gw)}\gd(x_*)^{2-\ga_-}), $$
where $c_3=c_3(N,\gm,c_1)$.
Therefore
$$ |\nabla \zeta(x)| \leq c_4\gd(x)^{\ga_+-1} \forevery x \in \Gw_{\gb_0} $$
where $c_4=c_4(N,\gm,c_1,\norm{f}_{L^\infty(\Gw)})$. Thus $\gd^{-\ga_-}|\nabla
\zeta| \in L^\infty(\Gw)$. \qed \medskip

\bdef{regularization}
Let $x_0 \in \Gw$ and denote $\tl\gb(x_0)=\min(\gb_0,\rec{2}\gd(x_0))$. We say
that $\tl G_\gm^\Gw$ is a \emph{proper regularization} of $G_\gm^\Gw$ relative to $x_0$ if $\tl G_\gm^\Gw(x)=G_\gm^\Gw(x_0,x)$ for $x \in \ovl \Gw_{\tl\gb(x_0)}$, $\tl
G_\gm^\Gw \in C^2(\Gw)\cap C(\ovl \Gw)$ and $\tl G_\gm^\Gw \geq 0$ in $\Gw$.
Similarly $\tl \gd$ is a \emph{proper regularization} of $\gd$ relative to $x_0$ if $\tl
\gd(x)=\gd(x)$ for $x \in \ovl \Gw_{\tl\gb(x_0)}$, $\tl \gd \in C^2(\ovl \Gw)$ and
$\tl \gd \geq 0$ in $\Gw$.
\es

\Remark Using \eqref{Gmu} and \eqref{vgf1}, it is easily  verified that
the functions $\vgf_{\mu,1}$, $\BBG_\gm^\Gw[\eta]$ (for $\eta \in
L^\infty(\Gw)$), $\tl G_\gm^\Gw$ and $\tl\gd^\gap$ belong to $X(\Gw)$.
Moreover, using \rprop{testgrad}, one obtains,
$$\gz\in X(\Gw)\q\text{and}\q h\in C^2(\bar\Gw)\Longrightarrow h\gz\in X(\Gw).$$
\vskip 3mm

In the proofs of the next two propositions we use the following construction. Let  $D\Subset\Gw$ be a $C^2$ domain.  The Green function for  $-L_\mu$ in $D$ is denoted by $G_\gm^{D}$. (To avoid misunderstanding we point out that, in the formula defining $L_\mu$, $\gd(x)$ denotes, as before, the distance from $x$ to $\bdw$, not to $\prt D$.) Given $x_0\in \Gw$  we construct a family of functions  $\CG(x_0)=\{\tl G_\gm^{D_\gb}:0<\gb<\rec{2}\tl\gb(x_0)\}$ such that, for each $\gb$, $\tl G_\gm^{D_\gb}$ is a proper regularization of $G_\gm^{D_\gb}(x_0,\cdot)$ in $D_\gb$ and $\CG(x_0)$ has the following properties:
\begin{itemize}
\item For every $\gb \in (0,\rec{2}\tl\gb(x_0))$, $\tl G_\gm^{D_\gb} \in C^2(\ovl D_\gb)$,
$\tl G_\gm^{D_\gb} \geq 0$ and
$\tl G_\gm^{D_\gb}(x)=G_\gm^{D_\gb}(x_0,x)$ for $x \in D_\gb\sms D_{\tl\gb(x_0)}$, .
\item The sequences $\{\tl G_\gm^{D_\gb}\}$ and $\{L_\gm\tl G_\gm^{D_\gb}\}$ converge to $\tl G_\gm^\Gw$ and $L_\gm \tl G_\gm^{\Gw}$ respectively, as $\gb \to 0$, a.e. in $\Gw$.
\item $\norm{\tl G_\gm^{D_\gb}+|L_\gm \tl G_\gm^{D_\gb}|}_{L^\infty(D_\gb)} \leq M_{x_0}$ where $M_{x_0}$ is a positive constant independent of $\gb$.
\end{itemize}
$\CG(x_0)$ will be called a \emph{uniform regularization} of $\{G_\mu^{D_\gb}\}$.

For any function $h \in C^2(\prt \Gw)$, we say that $\tl h$ is an \emph{admissible
extension} of $h$ relative to $x_0$ in $\ovl \Gw$ if $\tl h(x)=h(\gs(x))$ for
$x \in \Gw_{\tl\gb(x_0)}$ and $\tl h \in C^2(\ovl \Gw)$.

\subsection{Nonhomogeneous linear equations} Here we discuss the boundary value problem \eqref{eqG} in $\Gw$.

\blemma{intform} Let $u\in L^1_{loc}(\Gw)$ be a positive solution (in the sense of distributions) of equation
\bel{nonho} -L_\gm u = \gt \ee
in $\Gw$ where $\gt$ is a non-negative Radon measure.

If $\tau\in \GTM_{\gd^{\ga_+}}(\Gw)$ then
\begin{equation}\label{repG}
  -\int_{\Gw}\BBG_\gm^\Gw[\gt] L_\gm\zeta dx= \int_{\Gw} \zeta
d\gt \forevery \zeta \in X(\Gw).
\end{equation}
\es
\proof We may assume that $\tau$ is positive. By \rprop{Gmt}, $tr^*(\BBG_\gm^\Gw[\gt])=0$.
Therefore, given $\vge>0$, there exists $\bar \gb =\bar \gb(\vge)<\rec{2}\gb_0$  such that,
\begin{equation}\label{ep1}\BAL
\rec{\gb^\gan}\int_{\Gs_\gb}\BBG_\gm^\Gw[\gt]dS_x<\vge \q\text{and}\q
\int_{\Gw_{\gb}}\gd^{\ga_+}d\gt < \vge\q \forall \gb\in(0,\bar\gb].
\EAL\end{equation}
Let
$$I(\gb):=\int_{D_\gb} \BBG_\gm^\Gw[\gt]L_\gm \zeta dx +
\int_{D_\gb}\zeta d\gt.$$
To prove \eqref{repG} we show that
\begin{equation}\label{Igb}
\lim_{\gb\to0}I(\gb)=0.
\end{equation}

Put
$$\gt_1:=\chi_{_{\bar D_{\bar\gb}}}\gt,\q \gt_2:=\chi_{_{\Gw_{\bar\gb}}}\gt$$
and, for $0<\gb<\bar \gb$,
$$I_k(\gb):=\int_{D_\gb} \BBG_\gm^\Gw[\gt_k]L_\gm \zeta dx +
\int_{D_\gb}\zeta d\gt_k,\q k=1,2.$$
As $|\gz|\leq c\gd^\gap$ and $|L_\mu\gz|\leq \frac{c}{\gd^\gan}$, \eqref{ep1} implies,
\begin{equation}\label{I2gb}
 |I_2(\gb)|\leq c\vge \q\forall \gb\in (0,\bar\gb).
\end{equation}

For every $\gb \in (0,\bar\gb)$,
$$ -\int_{D_\gb}{}\BBG_\gm^\Gw[\gt_1]L_\gm \zeta dx = \int_{D_\gb}{}\zeta
d\gt_1 + \int_{\Gs_\gb}{}\frac{\prt \BBG_\gm^\Gw[\gt_1]}{\prt {\bf n}}
\zeta dS_x - \int_{\Gs_\gb}{} \BBG_\gm^\Gw[\gt_1]\myfrac{\prt \zeta}{\prt
\bf{n}}dS_x.
$$
Thus
$$I_1(\gb)=-\int_{\Gs_\gb}\frac{\prt \BBG_\gm^\Gw[\gt_1]}{\prt {\bf n}}
\zeta dS_x+ \int_{\Gs_\gb}{} \BBG_\gm^\Gw[\gt_1]\myfrac{\prt \zeta}{\prt
\bf{n}}dS_x=:I_{1,1}(\gb)+I_{1,2}(\gb).
$$
By \rprop{testgrad} and \eqref{ep1},
\begin{equation}\label{I12gb}
 |I_{1,2}(\gb)|\leq c\vge \q \forall \gb\in (0,\bar\gb).
\end{equation}

Next we estimate $I_{1,1}(\gb)$ for $\gb\in(0,\bar\gb/2)$. By Fubini,

$$ \BAL
I_{1,1}(\gb)&= -\int_{\Gs_\gb}\frac{\prt }{\prt \mathbf n_x} \int_{D_{\bar\gb}}G_\mu^\Gw(x,y)d\tau_1(y)\gz(x)dS_x\\
&= -\int_{D_{\bar\gb}}\int_{\Gs_\gb}\frac{\prt G_\mu^\Gw(x,y)}{\prt \mathbf n_x}
\gz(x)\,dS_x d\tau_1(y).
\EAL $$

For every $y\in D_{\bar\gb}$ the function $x\mapsto G_\mu^\Gw(x,y)$ is $L_\mu$-harmonic in $\Gw_{\bar\gb}$. By local elliptic estimates, for every $\gx\in \Gs_\gb$,
$$\sup_{x\in B_{\gb/4}(\gx)}|\nabla_x G_\mu^\Gw(x,y)|\leq c \gb^{-1} \sup_{x \in B_{\gb/2}(\gx)}G_\mu^\Gw(x,y).$$
By Harnack's inequality,
$$\sup_{x \in B_{\gb/2}(\gx)}G_\mu^\Gw(x,y)\leq c' \inf_{x \in B_{\gb/2}(\gx)}G_\mu^\Gw(x,y).$$
The constants $c,c'$ are independent of $\gb\in (0,\bar\gb/2)$, $y\in D_{\bar\gb}$ and $\gx\in \Gs_\gb$. Therefore we obtain,
\begin{equation}\label{gradest}
 |\nabla_x G_\mu^\Gw(x,y)|\leq C \gb^{-1}G_\mu^\Gw(x,y)\q \forall x\in \Gs_\gb,\;\forall y\in D_{\bar\gb},\;\forall \gb\in (0,\bar\gb/2).
\end{equation}
Hence,
$$|I_{1,1}(\gb)|\leq C\gb^{-1}\int_{\Gs_\gb}\BBG_\mu^\Gw[\tau_1]|\gz| \,dS_x.$$
As $|\gz(x)|\leq c\gd(x)^\gap$ it follows that,
$$|I_{1,1}(\gb)|\leq C\rec{\gb^{\gan}}\int_{\Gs_\gb}\BBG_\mu^\Gw[\tau_1]dS_x.$$
Therefore, by \eqref{ep1},
\begin{equation}\label{I11gb}
|I_{1,1}(\gb)|\leq C'\vge\q \forall \gb\in (0,\bar\gb/2).
\end{equation}
Finally \eqref{Igb} follows from \eqref{I2gb}, \eqref{I12gb} and \eqref{I11gb}.
 \qed

\bth{solG} Let $\gn \in \GTM^+(\prt \Gw)$ and $\gt \in
\GTM_{\gd^{\ga_+}}(\Gw)$. Then:

(i) Problem \eqref{eqG}has a unique solution.
The solution is given by
\bel{fo2} u=\BBG_\gm^\Gw[\gt] + \BBK_\gm^\Gw[\gn]. \ee

(ii) There exists a positive constant $c=c(N,\gm,\Gw)$ such that
\bel{estnorm} \norm{u}_{L^1_{\gd^{-\ga_-}}(\Gw)} \leq c(\norm{\gt}_{\GTM_{\gd^{\ga_+}}(\Gw)} +
\norm{\gn}_{\GTM(\prt \Gw)}). \ee

(iii) $u$ is a solution of of \eqref{eqG}   if and only
if $u \in L^1_{\gd^{-\ga_-}}(\Gw)$ and
\bel{repr} -\myint{\Gw}{}u L_\gm\zeta dx= \myint{\Gw}{} \zeta d\gt -
\myint{\Gw}{}\BBK_\gm^\Gw[\gn]L_\gm \zeta dx \forevery \zeta \in X(\Gw). \ee
\es

\proof (i) \rprop{Gmt} implies that \eqref{fo2} is a solution of \eqref{eqG}.

If $u$ and $u'$ are two solutions of \eqref{eqG} then $v:=(u-u')_+$ is a
nonnegative $L_\gm$-subharmonic function such that $tr^*(v)=0$ and $v\leq 2\BBG_\mu^\Gw[|\tau|]$ which is  a positive $L_\mu$-superharmonic function. By \rprop{subhar}, $v\equiv 0$ and hence $u \leq u'$ in $\Gw$. Similarly $u'\leq u$, so that
$u=u'$.
\medskip

(ii) In view of \eqref{estG1} and \eqref{estP}, \eqref{estnorm}  is an immediate consequence of \eqref{fo2}.
\medskip

(iii) Let $u$ be the solution of \eqref{eqG}. By \eqref{estnorm},
$u\in L^1_{\gd^{-\ga_-}}(\Gw)$ and by \rlemma{intform} and \eqref{fo2}, $u$ satisfies \eqref{repr}.

Conversely, suppose that $u \in L^1_{\gd^{-\ga_-}}(\Gw)$ and satisfies
\eqref{repr}. We show that $u$ is a solution of \eqref{eqG} or, equivalently, of \eqref{fo2}.

By \eqref{repr} with $\gz\in C_c^\infty(\Gw)$, $u$ is a solution (in the sense of distributions) of the equation in \eqref{eqG}. It remains to show that $tr^*(u)=\nu$.
 Put $U=u-\BBG_\gm^\Gw[\gt] - \BBK_\gm^\Gw[\gn]$ and note that, as $-L_\mu u=\tau$, $U$ is $L_\mu$-harmonic.

  Let $z \in \Gw$ and let $\CG(z)$ be a uniform regularization of $\{G_\mu^{D_\gb}:0<\gb<\rec{2}\tl\gb(z)\}$, (see Section \ref{S:test}). Then, for every $ \gb \in (0,\rec{2}\tl\gb(z))$, $\tl G_\mu^{D_\gb}\in C^2_0(\ovl D_\gb)$. Recall that $\tl G_\mu^{D_\gb}(x)=G_\mu^{D_\gb}(z,x)$.
Therefore, as $\frac{\prt G_\mu^{D_\gb}(z,x)}{\prt\mathbf n_x}=P_\mu^{D_\gb}(z,x)$, $x\in \Gs_\gb$, we obtain
\bel{w3} - \int_{D_\gb} U(x) L_\gm \tl G_\gm^{D_\gb}(x) dx =
\int_{\Gs_\gb}U(x) P_\gm^{D_\gb}(z,x)dS_x=U(z). \ee
The second equality is a consequence of the fact that $U$ is $L_\mu$-harmonic.
But $L_\gm \tl G_\gm^{D_\gb}(x)\to L_\gm \tl G_\mu^\Gw(z,x)$ pointwise and the sequence $\{L_\gm \tl G_\gm^{D_\gb}\}$ is bounded by a constant $M_z$. We observe that $U\in L^1(\Gw)$; in fact by assumption $u\in L^1_{\gd^{-\ga_-}}(\Gw)$ and therefore, by \rprop{GP}, $U\in L^1_{\gd^{-\ga_-}}(\Gw)$. Consequently, by \eqref{w3},
$$U(z)= - \int_{\Gw} U(x) L_\gm \tl G_\gm^{\Gw}(z,x) dx.$$
Since $G_\gm^{\Gw}(z,\cdot)\in X(\Gw)$, by \eqref{repr}  the right hand side vanishes. Thus $U$ vanishes in $\Gw$, i.e., $u$ satisfies \eqref{fo2}.
\qed

\bcor{supersol} \textit{Let $u$ be a positive $L_\mu$ superharmonic function. Then there exist
$\nu\in \GTM^+(\bdw)$ and $\tau\in \GTM^+_{\gd^\gap}(\Gw)$ such that \eqref{rep1} holds.
}
\es

\proof By the Riesz decomposition theorem $u$ can be written in the form $u=u_p+u_h$ where $u_p$ is an $L_\mu$-potential and $u_h$ is a non-negative $L_\mu$-harmonic function. Therefore there exists $\nu\in \GTM^+(\bdw)$ such  that   $u_h=\BBK_\mu^\Gw[\nu]$. Since $u_p$ is an $L_\mu$-potential there exists a positive Radon measure $\tau$ such that $u_p=\BBG_\mu^\Gw[\tau]$ (see e,g. \cite[Theorem 12]{An1}). This necessarily implies that $\tau\in \GTM_{\gd^\gap}(\Gw)$.
\qed

\bprop{subhar'}
Let $w$ be a non-negative $L_\mu$-subharmonic function. If $w$ has a normalized boundary trace then  it is dominated by an $L_\mu$-harmonic function.
\es
\proof
There exist a positive Radon measure $\tau$ in $\Gw$ and a measure $\nu\in \GTM^+(\bdw)$ such that
$$-L_\mu w=-\tau\q\text{in }\Gw,\q tr^*(w)=\nu.$$
 Let $u_\gb$ be the solution of
$$ -L_\mu u=-\tau_\gb\q \text{in }D_\gb,\q u=\BBK_\mu^\Gw[\nu]\q \text{on } \Gs_\gb$$
where $\tau_\gb:=\tau\chi_{_{D_\gb}}$. Then,
$$u_\gb+\BBG_\mu^{D_\gb}[\tau_\gb]=\BBK_\mu^\Gw[\nu].$$
Letting $\gb\to0$ we obtain,
$$\BBG_\mu^{\Gw}[\tau]\leq \BBK_\mu^\Gw[\nu].$$
Hence $\tau\in \GTM_{\gd^\gap}^+(\Gw)$ and consequently
\begin{equation}\label{subhar'}
 w+\BBG_\mu^{\Gw}[\tau]=\BBK_\mu^\Gw[\nu].
\end{equation}
\qed
\medskip

\section{The nonlinear equation}
In this section, we consider the nonlinear equation
\bel{nonl} - L_\gm u + u^q = 0 \ee
in $\Gw$ with $0<\gm<C_H(\Gw)$ and $ q>1$.  \medskip

{\bf Proof of Theorem A.} Since $u$ is a positive solution of \eqref{N}, $u$ is $L_\mu$-subharmonic. Assuming (i), $u$ is dominated by an $L_\mu$-harmonic function. Therefore, by \rprop{subhar}, $(i)\Lra (ii)$ and $u\in L^q_{\gd^\gap}(\Gw)$. On the other hand, by \rprop{subhar'} $(ii)\Lra (i)$.

As mentioned above, (i) implies that $u\in L^q_{\gd^\gap}(\Gw)$ and that there exists $\nu\in \GTM^+_{\gd^\gap}(\prt \Gw)$ such that $tr^* (u)=\nu$.  Therefore, by \rth{solG}, \eqref{form} is a consequence of \eqref{repr}. Thus $(i)\Lra(iii)$.

Finally, the implication $(iii)\Lra(i)$ is obvious.
\medskip

It remains to prove the last assertion. If $u$ is a positive solution of \eqref{PN} then, by (iii), $u\in L^q_{\gd^\gap}(\Gw)$ and \eqref{intrep} follows from \rth{solG}.

Conversely,  assume that $\gd^\gap u^q,\,u/\gd^\gan\in L^1(\Gw)$ and \eqref{intrep} holds. Then, by  \eqref{intrep} with $\gz\in C_c^\infty(\Gw)$, $u$ is a solution of \eqref{N}. Taking $\gz_f=\BBG_\mu^\Gw[f]$ where $f\in C_c(\Gw)$ and $f\geq 0$  we obtain
$$\int_\Gw (\BBK_\mu^\Gw[\nu]-u)f\;dx = \int_\Gw u^q\gz_f\;dx<\infty.$$
This implies $u\leq \BBK_\mu^\Gw[\nu]$, i.e., $u$ is $L_\mu$-moderate. Therefore by (i), $u$ is a solution of \eqref{PN}.
\qed

\bigskip

{\bf Proof of Theorem B.}

\note{Uniqueness} Let $u_1$ and $u_2$  be two positive
solutions of \eqref{PN}. Then $v:=(u_1-u_2)_+$ is a subsolution of \eqref{N} and therefore an $L_\mu$-subharmonic function. Furthermore, by (iii) in Theorem A, $u_1,u_2\in L^q_{\gd^\gap}(\Gw)$ and  $v\leq \BBG_\mu^\Gw[u_1^q+u_2^q]=:\bar v$. Obviously  $\bar v$ is $L_\mu$ superharmonic and $tr^*(v)=0$. Therefore, by \rprop{subhar},
$v=0$. Thus $u_1\leq u_2$ and similarly $u_2\leq u_1$.

 \medskip

\note{Monotonicity} As before, $v:=(u_1-u_2)_+$ is $L_\mu$-subharmonic and it is dominated by an $L_\mu$-superharmonic function. Since $\nu_1\leq \nu_2$, $tr^*(v)=0$. Hence by \rprop{subhar}, $v=0$.

\note{A-priori estimate} Suppose that $u$ is a positive solution of
\eqref{PN}.  Then \eqref{intrep} with $\zeta=\BBG_\gm^\Gw[1]$ implies \eqref{estnormN}. (Recall that $\BBG_\gm^\Gw[1] \sim \gd^{\ga_+}$.)
\qed
\medskip

For the proof of the next theorem we need

\blemma{ex_in_D} Let $D\Subset \Gw$ be a $C^2$ domain and $q>1$. If $h$ is a positive function in $L^1(\prt D)$ then there exists a unique solution of the boundary value problem,
\begin{equation}\label{bvpD}\BAL
 -L_\mu u + u^q&=0 \q \text{in }D\\
 u&=h\q \text{on }\prt D.
\EAL\end{equation}
\es

\proof First assume that $h$ is bounded. Let $P_{\mu}^D$ denote the Poisson kernel of $-L_\mu$ in $D$ and put $u_0:=\BBP_\mu^D[h]$. Thus $u_0$ is bounded. We show that there exists a non-increasing sequence of positive functions $\{u_n\}_1^\infty$, dominated by $u_0$, such that  $u_n$ is the solution of the boundary value problem,
\begin{equation}\label{bvpDn}\BAL
 -\Gd  v  + v^q&=\frac{\mu}{\gd^2}u_{n-1}\q \text{in }D\\
 v&=h\q \text{on }\prt D\q n=1,2,\ldots
\EAL\end{equation}
As usual $\gd$ denotes the distance to $\bdw$, not to $\prt D$. For $n=1$, $u_0$ is a supersolution of the problem and, obviously $v=0$ is a subsolution. Consequently there exists a  unique solution $u_1$. By induction, for $n>1$,
$$-\Gd  u_{n-1} + u_{n-1}^q=\frac{\mu}{\gd^2}u_{n-2}\geq \frac{\mu}{\gd^2}u_{n-1}.$$
Thus $v=u_{n-1}$ is a supersolution of \eqref{bvpDn} and it is bounded. It follows that there exists $0\leq u_n\leq u_{n-1}$ such that
$$-\Gd  u_{n} + u_{n}^q=\frac{\mu}{\gd^2}u_{n-1}\;\text{in }D,\q u_n=h\;\text{on }\prt D.$$
As the sequence is monotone we conclude that $u=\lim u_n$ is a solution of \eqref{bvpD}.

If $h\in L^1(\prt D)$, we approximate it by a monotone increasing sequence of non-negative bounded functions $\{h_k\}$. If $v_k$ is the solution of \eqref{bvpD} with $h$ replaced by $h_k$ then $\{v_k\}$ increases (by the comparison principle  \cite[Lemma 3.2]{BMR}) and $v=\lim v_k$ is a solution of \eqref{bvpD}.

Uniqueness follows by the comparison principle.
\qed

\medskip

{\bf Proof of Theorem C.} Put $u_0:=\BBK_\gm^\Gw[\gn]$ and $h_\gb:=u_0\lfloor_{\Gs_\gb}$.
Let $u_\gb$ be the solution of \eqref{bvpD} with $h$ replaced by $h_\gb$, $\gb\in (0,\gb_0)$. Since $u_0$ is a supersolution of \eqref{N} it follows that $\{u_\gb\}$ decreases as $\gb\downarrow 0$. Therefore $u:=\lim_{\gb\to 0} u_\gb$ is a solution of \eqref{N}.

We claim that $tr^*(u)=\nu$. Indeed,
\begin{equation}\label{C-gb}
 u_\gb+\BBG^{D_\gb}_\mu[u_\gb^q]=\BBP^{D_\gb}_\mu[h_\gb]=u_0.
\end{equation}

Furthermore, in $D_\gb$, $u_\gb\leq u_0\in L^q_{\gd^\gap}(\Gw)$. Therefore
$$\BBG^{D_\gb}_\mu[u_\gb^q]\to \BBG^\Gw_\mu[u^q].$$
Hence, by \eqref{C-gb},
$$u+\BBG^{\Gw}_\mu[u^q]=u_0=\BBK_\gm^\Gw[\gn].$$
By \rprop{Gmt}, $tr^*(u)=\nu$.

By Theorem B the solution is unique.

\qed
\medskip

\noindent{\bf Proof of Corollary C1.} By the previous theorem, if $\nu=f$ where $f$ is a positive bounded function then \eqref{PN} has a solution.  If $0\leq f\in L^1(\Gw)$ then it is the limit of an increasing sequence of such functions. Therefore, once again problem \eqref{PN} with $\nu=f$ has a solution.
\medskip

\noindent{\bf Proof of Theorem D.}\hskip 2mm Put $v=\BBK_\mu^\Gw[\gn]-u$. By the comparison principle $v\geq 0$. Clearly $v$ is $L_\gm$-superharmonic
in $\Gw$ and, by definition $tr^*(v)=0$. By \rth{Ri}, $v$ can be written under the form $ v_\gm= v_1
+ v_2 $ where $v_1$ is a nonnegative $L_\gm$-harmonic function and $v_2$ is an
$L_\gm$-potential. Since  $tr^*(v)=0$, it follows that
$tr^*(v_1)=tr^*(v_2)=0$. By the representation theorem, $v_1 = 0$. Therefore
$v=v_2$ and \eqref{behavior} follows by \rth{pot}. \qed

\bigskip

\noindent{\bf Proof of Theorem E.} By \rprop{GP}, specifically inequality \eqref{estP}, $\BBK_\mu^\Gw[\gn] \in L^q_{\gd^\gap}(\Gw)$ for every $q\in (1,q_{\mu,c})$ and $\nu\in \GTM^+(\bdw)$. Therefore the first assertion of the theorem is a consequence of Theorem C.

We turn to the proof of stability. Put $v_n=\BBK_\mu^\Gw[\nu_n]$. By \rprop{GP}, $\{v_n\}$ is bounded in $L^q_{\gd^\gap}(\Gw)$ for every $q\in (1,q_{\mu,c})$ and in $L^p_{\gd^{-\gan}}(\Gw)$ for every $p\in (1,\frac{N-\gan}{N-1-\gan})$. In addition $v_n\to v$ pointwise in $\Gw$.

This implies that $\{v_n^q\gd^\gap\}$ and $\{v_n/\gd^\gan\}$ are uniformly integrable in $\Gw$. Since $u_{\nu_n}\leq v_n$ it follows that this conclusion applies also to $\{u_{\nu_n}\}$.

 By the extension of the Keller -- Osserman inequality due to \cite{BMR}, the sequence $\{u_{\nu_n}\}$ is uniformly bounded in every compact subset of $\Gw$. Therefore, by a standard argument, we can extract a subsequence, still denoted by $\{u_{\nu_n}\}$ that converges pointwise to a solution $u$ of \eqref{N}.
In view of the uniform convergence mentioned above we conclude that
$$u_{\nu_n}\to u \q \text{in}\; L^q_{\gd^\gap}(\Gw) \text{ and in } L^1_{\gd^{-\gan}}(\Gw).$$
By Theorem A,
$$u_{\nu_n} + \BBG_\mu^\Gw[u_{\nu_n}^q]= \BBK_\mu^\Gw[\nu_n].$$
In view of the previous observations, passing to the limit as $n\to\infty$, we obtain,
$$u + \BBG_\mu^\Gw[u^q]= \BBK_\mu^\Gw[\nu].$$
Again by Theorem A it follows that $u$ is the (unique) solution of \eqref{PN}. Because of the uniqueness we conclude that the entire sequence $\{u_{\nu_n}\}$ (not just a subsequence) converges to $u$ as stated in assertion II. of the theorem.

Finally we prove assertion III. By Theorem A
\bel{formula} u_{k\gd_y} + \BBG_\gm^\Gw[u_{k\gd_y}^q]=kK_\gm^\Gw(\cdot,y). \ee
Combining \eqref{Pmu}, \eqref{Gmu} and the fact $u_{k\gd_y} \leq
kP_\gm^\Gw(\cdot,y)$, we obtain
$$ \frac{\BBG_\gm^\Gw[u_{k\gd_y}^q](x)}{K_\gm^\Gw(x,y)} \leq k^q
\frac{\BBG_\gm^\Gw[(K_\gm^\Gw(.,y)^q](x)}{K_\gm^\Gw(x,y)} \leq
ck^q|x-y|^{N+\ga_+-q(N-1-\ga_-)}.
$$
Since $1<q<q_{\gm,c}$, it follows that
$$ \lim_{x \to 0}\myfrac{\BBG_\gm^\Gw[u_{k\gd_y}^q](x)}{K_\gm^\Gw(x,0)}=0. $$
Therefore, by \eqref{formula}, we obtain \eqref{dirac}. \qed \medskip
\medskip

\noindent{\bf Proof of Theorem F.} Let $y \in \prt \Gw$. By negation, assume that there exists a positive solution $u$ of
\eqref{PN} with $\nu=k\gd_y$ for some $k>0$.
By Theorem A, $u \leq k\BBK_\gm^{\Gw}(.,y)$ and $u \in L^q_{\gd^{\ga_+}}(\Gw)$.
Let $\gg\in (0,1)$ and denote $C_\gg(y)=\{x \in \Gw: \gg|x-y| \leq \gd(x)\}$. By Theorem E III.,
$$ \lim_{x \in C_\gg(y), x\to y}\myfrac{u(x)}{K_\gm^\Gw(x,y)}=k. $$
This implies that there
exist positive numbers $r_0, c$  such that
\bel{contr} u(x) \geq c K_\gm^\Gw(x,y) \forevery x \in C_\gg(y) \cap
B_{r_0}(y).\ee
By \eqref{Pmu},
$$ \BA{lll}J_\gg:= \int_{C_\gg(y) \cap B_{r_0}(y)}(K_\gm^\Gw(x,y))^q\gd(x)^{\ga_+}dx \\ [3mm]
\phantom{------}
\geq c'\int_{C_\gg(y) \cap B_{r_0}(y)}{}\gd(x)^{\ga_+(q+1)}|x-y|^{(2\ga_- -N)q}dx \\ [3mm]
\phantom{------}
\geq c'\gg^{\ga_+(q+1)}\int_{C_\gg(y) \cap B_{r_0}(y)}|x-y|^{\ga_+ -
q(N-1-\ga_-)}dx.
\EA $$
Since $q \geq q_{\gm,c}$ the last integral is divergent.  But \eqref{contr} and the fact that $u \in L^q_{\gd^{\ga_+}}(\Gw)$ imply that $J_\gg<\infty$.
We reached a contradiction. \qed \bigskip

\noindent{\bf Acknowledgements} This research was supported by the Israel
Science Foundation founded by the Israel Academy of Sciences and Humanities,
through grant 91/10. Part of this research was carried out, by the first author, during a visit at the Isaac Newton Institute, Cambridge as part of the FRB program. He wishes to thank the institute for providing a pleasant and stimulating atmosphere.  The second author was also partially supported  by a Technion
fellowship.

The authors wish to thank Professor Pinchover for many useful discussions.

\end{document}